\documentclass[12pt]{amsart}
\usepackage{amsmath,amsthm,amsfonts,amssymb}
\usepackage[all]{xy}
\usepackage{amssymb, amsthm, amsmath}
\usepackage[english]{babel}
\usepackage{amsfonts}
\newcommand{\ra}{\rightarrow}

\newcommand{\ot}{\otimes}
\newcommand{\mtc}{\mathcal}
\newcommand{\lam}{\lambda}
\newcommand{\Lam}{\Lambda}

\newcommand{\al}{\alpha}
\newcommand{\eps}{\epsilon}
\newcommand{\bn}{\begin}
\newcommand{\en}{\end}

\newcommand{\sub}{\subsection}

\newcommand{\D}{\Delta}

\numberwithin{equation}{section}

\newtheorem{thm}[equation]{Theorem}
\newtheorem{prop}[equation]{Proposition}

\newtheorem{cor}[equation]{Corollary}
\newtheorem{rem}[equation]{Remark}

\title[Normal Hopf subalgebras]
{Normal Hopf Subalgebras of Semisimple Hopf Algebras}

\newcommand{\AMSclasifR}{16W35, \;16W40}

\author{Sebastian  Burciu}
\address{Inst.\ of Math.\ ``Simion Stoilow" of the Romanian Academy
P.O. Box 1-764, RO-014700, Bucharest, Romania\\ smburciu@syr.edu}

\date{September 17, 2007}
\begin{document}
\thanks{MSC (2000): \AMSclasifR}\thanks{The research was supported by CEx05-D11-11/04.10.05.}
\begin{abstract} The notion of kernel of a representation of a semisimple Hopf algebra is introduced.
Similar properties to the kernel of a group representation are
proved in some special cases. In particular, every normal Hopf
subalgebra of a semisimple Hopf algebra $H$ is the kernel of a
representation of $H$. The maximal normal Hopf subalgebras of $H$ are described.
\end{abstract}
\maketitle

\newcommand{\lbd}{\Lambda}
\newcommand{\dw}{\downarrow}
\newcommand{\uw}{\uparrow}
\newcommand{\dlt}{\delta}
\newcommand{\nono}{\nonumber}
\newcommand{\ch}{\chi}
\newcommand{\mtr}{\mathrm}
\numberwithin{equation}{section}

\section*{Introduction}
In this paper the notion of kernel of a representation of a finite
dimensional semisimple Hopf algebra is proposed. Similar
properties to the kernel of a group representation are proved in
some special cases.

Let $G$ be a finite group and $ \mathrm{X}:G\rightarrow
\mathrm{End}_k(M)$ be a finite dimensional representation of $G$
which affords the character $\chi$. The kernel of the representation
$M$ is defined as $\mathrm{ker}\;\chi=\{g \in G|\;\ch(g)=\ch(1) \}$
and it is the set of all group elements $g \in G$ which act as
identity on $M$. (for example, see \cite{Is}) Every normal subgroup
$N$ of $G$ is the kernel of a character, namely the character of the
regular representation of $G/N$. If $\mathrm{Z}=\{g \in
G|\;\;|\ch(g)|=\ch(1)\}$ then $\mtr{Z}$ is called the $center$ of the
character $\chi$ and it is the set of group elements
of $G$ which act as a unit scalar on $M$. The properties of
$\mathrm{Z}$ and $\mathrm{ker}\;\chi$ are described in [Lemma 2.27,
\cite {Is}] which shows that $\mtr{Z}/ \mathrm{ker}\;\chi$ is a
cyclic subgroup of the center of $G/\mtr{ker}\;\ch$.


If $M$ is a representations of a finite dimensional semisimple Hopf algebra $H$ and $\chi \in C(H)$ is its associated character 
then $\mtr{ker}\;\chi \subset H$ is defined as the set of all irreducible
$H^*$-characters  $d \in H$ such that $d$ acts as the scalar
$\eps(d)$ on $M$. Note that in the case $H=kG$ the irreducible $H^*$-characters are the elements of $G$ and we obtain exactly the kernel
of $M$ as defined in \cite{Is}. We prove that $\mtr{ker}\;\chi=\{d
\in \mtr{Irr}(H^*)| \;\ch(d)=\eps(d)\ch(1)\}$. Similarly, the set of
all the irreducible characters of $H^*$ that acts as a root of unity scalar
on $M$ is characterized as $\mathrm{z}_{_\ch}=\{d \in
\mtr{Irr}(H^*)|\;\;|\ch(d)|=\eps(d)\ch(1)\}$.

Section \ref{propr} presents the definition and the main properties of the $kernel\;\ch$ of a character $\ch$ and its center $z_{_{\ch}}$. It is shown that these sets of characters are closed under multiplication and the duality operation $``\;^*\;"$ and thus they generate Hopf subalgebras of $H$ denoted by $H_{_{\ch}}$ and $Z_{_{\ch}}$, respectively. We say that a Hopf subalgebra $K$ of $H$ is the kernel of a representation if $K=H_{_{\ch}}$ for a certain character $\ch$ of $H$.

Section \ref{normal} studies the relationship between normal Hopf subalgebras and the Hopf algebras generated by kernels. It is shown that any normal Hopf subalgebra is the kernel of a character which is central in $H^*$.


Section \ref{centralch} investigates the structure of the Hopf subalgebras $H_{_{\ch}}$ generated by the kernel of a  character $\ch$ affording a representation $M$ and which is central in $H^*$. It is shown that $H_{_{\ch}}$ is normal in $H$ and all the representations of $H//H_{_{\ch}}:=H/HH_{_{\ch}}^+$ are the representations which are constituents of some tensor power of $M$. Combining this with the main theorem of the previous section it follows that a Hopf subalgebra is normal if and only if it is the kernel of a character central in the dual Hopf algebra.

Using a basis description given in \cite{Zh} for the algebra generated by the characters which are central characters in $H^*$  we describe a finite collection of normal Hopf subalgebras of $H$ which are the maximal normal Hopf subalgebras of $H$ ( under inclusion). Any other normal Hopf subalgebra is an intersection of some of these Hopf algebras. Two other results that also hold for group representation are presented in this section. In particular $Z_{_{\ch}}//H_{_{\ch}}$ is a group algebra of a cyclic group contained in the group of central grouplike elements of $H//H_{_{\ch}}$.

In the last section we consider a double coset type decomposition
for Hopf algebras. If $H$ is a finite dimensional semisimple Hopf algebra and $K$ and $L$ are two Hopf algebras then $H$ can be decomposed as sum of $K-L$ bimodules which are free both as $K$-modules and $L$-modules. To the end we give an application in the situation of a
unique double coset.

Algebras and coalgebras are defined over the algebraically closed ground field $k=\mathbb{C}$. For a vector space $V$ over $k$  by $|V|$ is denoted the dimension $\mtr{dim}_kV$. The comultiplication, counit and antipode of a
Hopf algebra are  denoted by $\Delta$, $\epsilon$ and
$S$, respectively. We use Sweedler's notation $\D(x)=\sum x_1\ot x_2$ for all $x\in H$.  All the other Hopf notations are those used in
\cite{Montg}.

\section{Properties of the kernel}\label{propr}

Let $H$ be a finite dimensional semisimple Hopf algebra over the algebraically closed field $k$ of characteristic zero. Then $H$ is also cosemisimple \cite{Lard}. Denote by $\mtr{Irr}(H)$ the set of irreducible characters of $H$. Let $\mathbb{G}_0(H)$ be the Grothendieck group of the category $H$-mod of finite dimensional left $H$-modules and
$C(H)$ be its scalar extension to $k$. Then $C(H)=\mathbb{G}_0(H)\ot_{_{\mathbb{Z}}} k$ is a semisimple subalgebra of $H^*$ with basis given by the characters of the irreducible $H$-modules \cite{Z}. Moreover $C(H)=\mtr{Cocom}(H^*)$, the space of cocommutative elements of $H^*$. The character ring of $H^*$ is a semisimple subalgebra of $H^{**}=H$ and under this identification it follows that $C(H^*)=\mtr{Cocom}(H)$.


If $M$ is an $H$-module with character $\chi$ then $M^*$ is also an $H$-module with character $\chi^*=\chi \circ S$. This induces an involution $``\;^*\;":C(H)\ra C(H)$ on $C(H)$.

Let $H$ be a finite dimensional semisimple Hopf algebra over an algebraically closed
field of characteristic zero. Recall that the exponent of $H$ is the
smallest positive number $m >0$ such that $h^{[m]}=\eps(h)1$ for all
$h \in H$. The generalized power $h^{[m]}$ is defined by
$h^{[m]}=\sum_{(h)}h_1h_2...h_m$. The exponent of a finite
dimensional semisimple Hopf algebra is always finite and divides the
cube power of dimension of $H$ \cite{EG'}.


If $W \in H^*$-mod then $W$ becomes a right $H$-comodule via $\rho :W \ra W \ot H$ given by $\rho(w)=\sum w_{0}\ot w_1$ if and only if $fw=\sum f(w_1)w_0$ for all $w \in W$ and $f \in H^*$.

\sub{} \label{assoc} Let $W$ be a simple $H^*$-module. Then $W$ is a simple right $H$-comodule and  one can associate to it a simple subcoalgebra of $H$ denoted by $C_{_W}$ \cite{Lar}. If $q=|W|$ then $|C_{_W}|=q^2$ and it is a co-matrix coalgebra. It has a basis $\{x_{ij}\}_{1\leq i,j \leq q}$  such that $\D(x_{ij})=\sum_{l=0}^qx_{il}\ot x_{lj}$ for all $1\leq i,j \leq q$. Moreover $W \cong k<x_{1i}|\;1\leq i \leq q>$ as right $H$-comodule where $\rho(x_{1i})=\D(x_{1i})=\sum_{l=0}^qx_{1l}\ot x_{li}$ for all $1\leq i \leq q$. The character of $W$ as left $H^*$-module is $d \in C(H^*)\subset H$ and it is given by $d=\sum_{i=1}^qx_{ii}$. Then $\eps(d)=q$ and the simple subcoalgebra $C_{ _W}$ is also denoted by $C_{ _d}$.
\begin{prop}\label{centraldacts} Let $H$ be a finite dimensional semisimple Hopf algebra over $k$ and $M$ and $W$ be irreducible representations of $H$ and $H^*$ respectively, affording the characters $\chi \in
C(H)$ and $d \in C(H^*)$. Then:
\begin{enumerate}\item $|\chi(d)| \leq
\chi(1)\eps(d)$ \item Equality holds if and only if $d$ acts as $\alpha
\eps(d) Id_{_M}$ on the irreducible \\$H$-representation $M$, where $\alpha$ is a root of
unity.
\end{enumerate}
\end{prop}

\begin{proof}
\begin{enumerate}
\item
$W$ is a right $H-$comodule and one can define the map
\begin{center}
\begin{tabular}{rclc}
$\mtr{T} :$ & $M \ot W$ &$\longrightarrow$ & $M \ot W$ \\
 & $m \ot w$ & $\longmapsto$ & $\sum w_1 m \ot w_0$\\
\end{tabular}
\end{center}
It can be checked that $T^p(m\ot w)=w_1^{[p]}m\ot w_0$ for all $p \geq 0$.
Thus, if $m=\mtr{exp}(H)$ then $\mtr{T}^m=\mtr{Id}_{_{M \ot W}}$. Therefore
$\mtr{T}$ is a semisimple operator and all its eigenvalues are root
of unity. It follows that $\mtr{tr}(\mtr{T})$ is the sum of all these
eigenvalues and in consequence $|\mtr{tr}(\mtr{T})| \leq \mtr{dim}_k
(M \ot W)= \chi(1)\eps(d)$.

It is easy to see that $\mtr{tr}(\mtr{T})=\chi(d)$. Indeed using the above remark one can suppose that $W=k<x_{1i}\;|\;1\leq i \leq q>$ where
$C_{_W}=k<x_{ij}|\; 1 \leq i,j \leq q >$ is the coalgebra associated to $W$. Then the formula for $T$ becomes
$\mtr{T}(m \ot x_{1i})=\sum_{j=0}^{\eps(d)} x_{ji}m \ot x_{1j}$ which shows that
$\mtr{tr}(\mtr{T})=\sum_{i=0}^{\eps(d)}\chi(x_{ii})=\chi(d)$.
\item

Equality holds if and only if $\mtr{T}=\alpha
\chi(1)\eps(d)\mtr{Id}_{_{M \ot W}}$ for some $\alpha$ root of
unity. The above expression for $\mtr{T}$ implies that in this case
$x_{ij}m=\delta_{i,j}\alpha m$ for any $1 \leq i,j \leq \eps(d)$. In
particular $dm=\alpha \eps(d) m$ for any $m \in M$ which shows that
$d$ acts as a scalar multiple on $M$ and that scalar is
$\alpha\eps(d)$. The converse is immediate.
\end{enumerate}
\end{proof}

Let $M$ be an irreducible representations of $H$ which affords the character
$\chi$. Define $\mtr{ker}\;\ch$ as the set of all irreducible characters $d \in \mtr{Irr}(H^*)$
which act as the scalar $\eps(d)$ on $M$. The previous proposition implies that $\mtr{ker} \;\chi=\{d \in \mtr{Irr}(H^*) |
\; \chi(d)=\eps(d)\chi(1)\}$.

Similarly let $\mtr{z}_{ _\ch}$ be the the set of all irreducible characters $d \in \mtr{Irr}(H^*)$ which act as a scalar
$\alpha\eps(d)$ on $M$, where where $\alpha$ is a root of unity. Then from the same proposition it follows
$\mtr{z}_{_\chi}=\{d \in \mtr{Irr}(H^*)
|\;\;|\chi(d)|=\eps(d)\chi(1)\}$.
Clearly $ \mtr{ker}\;\chi \subset \mtr{z}_{_{\chi}}$.

\begin{rem} \label{prod}For later use let us notice that $\mtr{ker}\;\chi \subset \mtr{ker}\;\chi^n$ for all $n \geq 0$. Indeed if $d\in C(H^*)$ is an element of $\mtr{ker}\;\chi$ then one has a simple subcoalgebra $C_d$ associated to $d$ and $d=\sum_{i=1}^{\eps(d)}x_{ii}$. Then from the proof of Proposition \ref{centraldacts} it follows that $\chi(x_{ij})=\chi(1)\delta_{ij}$. Thus $$\chi^n(d)=\sum_{i=1}^{\eps(d)}\sum_{i_1,\cdots,i_n=1}^{\eps(d)}\chi(x_{ii_1})\chi(x_{i_1i_2})\cdots \chi(x_{i_ni})=\chi(1)^n\eps(d)$$
Similarly it can be shown that $z_{_{\ch}}\subset z_{_{\ch^n}}$.
\end{rem}

\begin{rem}\label{frompf}
From the proof of Proposition \ref{centraldacts} if $d \in z_{_{\ch}}$ and $C_d=k<x_{ij}\; |\;1\leq i,j \leq \eps(d)>$ is the simple subcoalgebra of $H$ associated to $d$ then $\chi(x_{ij})=\chi(1)\al_{\ch}\delta_{ij}$ where $\al_{\ch} \in k$ is a root of unity.
\end{rem}


A subset $X \subset\mtr{Irr}(H^*)$ is closed under multiplication if for every $\chi, \mu \in X$ in the decomposition of $\chi\mu=\sum_{\gamma \in \mtr{Irr}(H)}m_{\gamma}\gamma$ one has $\gamma \in X$ if $m_{\gamma} \neq 0$. A subset $X \subset C(H)$ is closed under $``\;^*\;"$ if $x^* \in X$ for all $x \in X$.

\begin{prop}\label{clsets}
Let $H$ be a finite dimensional semisimple Hopf algebra and $M$ a simple $H$-module affording the irreducible character $\chi \in C(H)$. Then the subsets $\mtr{ker}\;\chi$ and $\mtr{z}_{_\chi}$ of $C(H^*)$ are
closed under multiplication and $``\;^*\;"$.
\end{prop}

\begin{proof}
Proposition \ref{centraldacts} implies that $\chi(d)=\eps(d)\chi(1)$
if and only if $d$ acts as $\eps(d)\mtr{Id}_{ _M}$ on $M$. Therefore if $d
\in \mtr{ker}\;\chi$ then $d^*=S(d) \in \mtr{ker} \;\chi$ since
$\chi(d^*)=\overline{\chi(d)}$ \cite{NR'}.

Let $d, d' \in \mtr{ker}\;\chi$. Then $dd'$ acts as $\eps(dd')Id_M$ on $M$ since $d$ acts as $\eps(d)Id_M$ and $d'$ acts as $\eps(d')Id_M$ on $M$
Write $dd'=\sum_{i=1}^qm_id_i$ where $d_i$ are irreducible
characters of $H^*$ and $m_i \neq 0$ for all $1\leq i \leq q$. Then $\chi(dd')=\sum_{i=1}^qm_i\chi(d_i)$ and
$$\chi(1)\eps(dd')=|\chi(dd')| \leq \sum_{i=1}^qm_i|\chi(d_i)| \leq
\chi(1)\sum_{i=1}^qm_i\eps(d_i)=\chi(1)\eps(dd')$$ It follows by
Proposition \ref{centraldacts} that $\chi(d_i)=\chi(1)\eps(d_i)$  and therefore $d_i
\in \mtr{ker}\;\chi$ for
all $1\leq i\leq q$.

Similarly it can be shown that $z_{ _{\chi}}$ is
closed under multiplication and " $^*$ ".
\end{proof}

\sub{} \label{hopfgen} If $X \subset C(H^*)$ is closed under multiplication and $``\;^*\;"$ then it generates a Hopf subalgebra of $H$ denoted by $H_{ _X}$ \cite{NR}. One has $H_{ _X}=\oplus_{d\in X}C_d$.
 Using this, since the sets $\mtr{ker}\;\chi$ and $\mtr{z}_{_\chi}$ are
closed under multiplication and $``\;^*\;"$ they generate Hopf
subalgebras $H$ denoted by $H_{_{\chi}}$ and  $\mtr{Z}_{_\chi}$, respectively.

\begin{rem} \label{restr} The proof of the Proposition \ref{centraldacts} implies that
$\chi\downarrow_{_{H_{_\chi}}}=\chi(1)\eps_{_{H_{\chi}}}$ where $\chi\downarrow_{_{H_{\chi}}}$ is the restriction of $\chi$ to the subalgebra $H_{_{\chi}}$ and $\eps_{_{H_{_{\chi}}}}$ is the character of the trivial module over the Hopf algebra $H_{_{\chi}}$.
\end{rem}

\section{Normal Hopf subalgebras}\label{normal}
As before let $H$ be a finite dimensional semisimple Hopf algebra over the field $k=\mathbb{C}$. In this section we prove that any normal Hopf subalgebra of $H$ is the kernel of a character which is central in the dual algebra $H^*$.
 
Let $\mtr{Irr}(H)=\{\ch_0,\cdots, \; \ch_s\}$ be the set of all irreducible $H$-modules. Let $M$ be  a not necessarily irreducible representations of $H$ which affords the character $\chi$. Define as before $$\mtr{ker} \;\chi=\{d \in \mtr{Irr}(H^*) | \; \chi(d)=\eps(d)\chi(1)\}$$ If $\chi=\sum_{i=1}^{p}m_i\chi_i$ where $m_i \in \mathbb{Z}_{\geq 0}$ then $\mtr{ker}\;\chi=\cap_{m_i\neq 0}\; \mtr{ker} \;\chi_i$, where $\mtr{ker} \; \chi_i$ was defined in the previous section.

Proposition \ref{centraldacts} implies that
$\mtr{ker}\;\ch$ is the set of all irreducible characters $d$ of
$H^*$ which act as the scalar $\eps(d)$ on $M$. If
$\mtr{z}_{_\chi}=\{d \in \mtr{Irr}(H^*)
|\;\;|\chi(d)|=\eps(d)\chi(1)\}$ then from the same proposition it follows $\mtr{z}_{_{\chi}}$ is the set of all
irreducible characters $d$ of $H^*$ which act as a scalar
$\alpha\eps(d)$ on $M$, where where $\alpha$ is a root of unity.
Clearly $ \mtr{ker}\;\chi \subset \mtr{z}_{_{\chi}} \subset \cap_{m_i\neq 0} \;z_{ _{\chi_i}} $. Proposition \ref{clsets} implies that $\mtr{ker}\;\ch$ and $\mtr{z}_{_{\chi}}$ are closed under multiplication and $``\;^*\;"$, thus they generate Hopf subalgebras of $H$ denoted again by $H_{ _{\ch}}$ and $\mtr{Z}_{ _\ch}$.

For a semisimple Hopf algebra $H$ over an algebraically closed field $k$ of characteristic zero use the notation $\Lam_H \in H$ for the idempotent integral of $H$ ( $\eps(\Lam_H)=1$) and $t_H \in H^*$ for the idempotent integral of $H^*$ ($t_H(1)=1)$.  If $\mtr{Irr}(H)=\{\ch_0,\cdots, \; \ch_s\}$ is the set of irreducible $H$-modules then from \cite{Montg} it follows that the regular character of $H$ is given by the formula \begin{equation}\label{f1}|H|t_{ _H}=\sum_{i=0}^s\chi_i(1)\chi_i\end{equation} The dual formula is
\begin{equation}\label{f2}|H|\Lambda_{ _H}=\sum_{d \in \mtr{Irr}(H^*)}\eps(d)d\end{equation}

One also has $t_{ _H}(\Lam_{ _H})=\frac{1}{|H|}$ \cite{Lard}.

If $K$ is a Hopf subalgebra of $H$ then $K$ is a semisimple and cosemisimple Hopf algebra \cite{Montg}.

A Hopf subalgebra $K$ of $H$ is called normal if $h_1xS(h_2)\in K$ and $ S(h_1)xh_2 \in K$ for all $x \in K$ and $h \in H$. If $H$ is semisimple Hopf algebra as above then $S^2=\mtr{Id}$ (see \cite{Lard}) and $K$ is normal in $H$ if and only if $h_1xS(h_2)\in K$ for all $x \in K$ and $h \in H$. Also, in this situation  $K$ is normal in $H$ if and only if $\Lam_{ _K}$ is central in $H$ \cite{Mas'}. If $K^+=Ker(\eps)\cap K$ and $K$ normal Hopf subalgebra of $H$ then $HK^+=K^+H$ and $H//K:=H/HK^+$ is a quotient Hopf algebra of $H$ via the canonical map $\pi:H\ra H//K$.

\sub{} \label{normalfacts} Suppose that $K$ is a normal Hopf subalgebra of $H$ and let $L=H//K$ be the quotient Hopf algebra of $H$ via $\pi:H\ra L$. Then $\pi^* :L^* \ra H^*$ is an injective Hopf algebra map. It follows that $\pi^*(L^*)$ is normal in $H^*$. Indeed, it can be checked that $\pi^*(L^*)=\{f\in H^*|f(ha)=f(h)\eps(a)\;\; \text{for all}\; h\in H,\;a\in K\}$. Using this description it is easy to see that $g_1fS(g_2)\in L^*$ for all $g \in H^*$ and $f\in L^*$. Moreover ${(H^*//L^*)}^*\cong K$ since ${(H^*//L^*)}^*=\{a\in H^{**}=H\;|\; fg(a)=f(1)g(a) \;\text{for all}\; f\in H^* \;\text{and}\; g \in L^*\}$.

In this situation it can be seen that the irreducible representations of $L=H//K$ are those representations $M$ of $H$ such that each $x\in K$ acts as $\eps(x)Id_M$ on $M$. Let $M$ be an irreducible representation of $L$. Suppose that $\chi$ is the character of $M$ as $L$-module. Then $\pi^*(\ch) \in C(H)$ is the character of $M$ as $H$-module. If $d \in C(K^*)\subset C(H^*)$ then $d$ acts as $\eps(d)\mtr{Id}_M$ on $M$ and therefore $d\in \mtr{ker}\; \chi$. Using the notations from the previous section it follows that $H_{_{\pi^*(\chi)}}\supset K$.

Let $\mu$ be any irreducible representation of $H$ and $\xi_{_\mu} \in \mtr{Z}(H)$ be the central primitive idempotent associated to it. Then $\nu(\xi_{\mu})=\delta_{\mu,\;\nu}\mu(1)$ for any other irreducible character of $H$.

Dually, since $H^*$ is semisimple to any irreducible character $d \in C(H^*)$ one has an associated central primitive idempotent $\xi_d \in H^*$. As before one can view $d \in H^{**}=H$ and the above relation becomes $\xi_d(d')=\delta_{d,\;d'}\eps(d)$ for any other irreducible character $d' \in C(H^*)$.

We say that a Hopf subalgebra $K$ of $H$ is the kernel of a character if $K=H_{_{\ch}}$ for some character $\ch \in C(H)$.
Following is the main result of this section.

\begin{thm}\label{main}
Let $H$ be a finite dimensional semisimple Hopf algebra over $k$. Any normal Hopf subalgebra of $H$ is the kernel of a character which is central
in $H^*$.
\end{thm}
\begin{proof}
Let $K$ be a normal Hopf subalgebra of $H$ and $L=H//K$. Then $L$ is a semisimple and cosemisimple Hopf algebra \cite{Montg}. The above remark shows that the irreducible representations of $L$ are exactly those irreducible
representations $M$ of $H$ such that $H_{_{\ch_{_M}}}
\supset K$ where $\ch_{_M}$ is the $H$-character of $M$.

Let $\pi :H \ra L$ be the natural projection and $\pi^* :L^* \ra H^*$ its dual map. Then $\pi^*$ is an injective Hopf algebra map and $L^*$ can be identified with a Hopf subalgebra of $H^*$.
Therefore, if $t_{_L} \in L^*$ is the idempotent
integral of $L$ then $|L|t_{ _L}$ is the regular representation of $L$ and $H_{|L|\pi^*( t_{ _L})} \supset K$. Since $\pi^*(L^*)$ is a normal Hopf subalgebra
of $H^*$ it follows from the same result \cite{Mas'} that $\pi^*(t_{_L})$ is a central element of $H^*$. Therefore if we show that $H_{_{|L|\pi^*(t_{_L})}} = K$ then the proof will be complete.

For any irreducible character $d \in C(H^*)$ let $\xi_d \in H^*$ be the associated central primitive of $H^*$.
Then $\{\xi_d\}_{d\in \mtr{Irr}(H^*)}$ is the complete set of central orthogonal idempotents of $H^*$ and $\xi_d(d')=\delta_{d,\;d'}\eps(d)$ for any two irreducible characters of $H^*$.

Since $\pi^*(t_{ _L})$ is a central idempotent of $H^*$ one can write it as a sum of central primitive orthogonal idempotents $$\pi^*(t_{_L})=\sum_{d \in X}\xi_d$$ where $X$ is a subset
of  $\mtr{Irr}(H^*)$. It follows that for any $d \in \mtr{Irr}(H^*)$ one has that $\pi^*(t_{ _L})(d)=\eps(d)$ if $d \in X$ and $\pi^*(t_{ _L})(d)=0$ otherwise.

Since $\mtr{ker}\;|L|\pi^*(t_{ _L}) \supset K$ one has $\pi^*(t_{_L})(d)=\eps(d)$ for all $d
\in \mtr{Irr}(K^*)$ and thus all the irreducible characters
of $\mtr{Irr}(K^*) \subset X$.

Let $\Lambda_{ _H} \in H$ and $\Lambda_{_L} \in L$ be
the idempotent integrals of $H$ and  $L$. Since $\pi$ is a Hopf algebra map one has $\pi(\Lam_H)=\Lam_{ _L}$.

Then $$\pi^*(t_{_L})(\Lambda_{ _H})=t_{ _L}(\pi(\Lambda_{ _H}))=t_{ _L}(\Lam_{ _L})=\frac{1}{|L|}$$

On the other hand,  since
$\Lambda_{ _H}= \frac{1}{|H|}\sum_{d \in \mtr{Irr}(H^*)}\eps(d)d$ it follows that
\begin{equation*}
\pi^*(t_{_L})(\Lambda_{ _H})=\frac{1}{|H|}\sum_{d
\in X}\eps^2(d)
\end{equation*}
which implies that $\sum_{d \in X}\eps^2(d)=\frac{|H|}{|L|}=|K|$.
Since $\sum_{d \in \mtr{Irr}(K^*)}\eps^2(d)=|K|$ and $\mtr{Irr}(K^*) \subset
X$ we conclude that $\mtr{Irr}(K^*)=X$ and $H_{_{|L|\pi^*(t_{_L})}}=K$
\end{proof}
\begin{rem}Let $A$ be a finite dimensional semisimple Hopf algebra over the field $k$ and $e \in A$ be an idempotent. Then $e$ is central in $A$ if and only if the value $\chi(e)$ is either $\chi(1)$ or $0$ for any irreducible character $\chi$ of $A$.

Indeed, since $k$ is algebraically closed $A$ is a product of matrix rings. If $e$ is a central idempotent then $e$ is a sum of primitive central idempotents and one has the above relations. Conversely, write $e$ as sum of primitive orthogonal idempotents of $A$. The above relations implies that $e$ contains either all primitive idempotents from a matrix ring or none of them. Thus $e$ is central in $A$.
\end{rem}

If $N$ is a $K$-module let $N\uparrow_K^H=H\ot_KN$ be the induced module. If $N$ has character $\mu \in C(K)$ then denote by $\mu\uparrow_K^H \in C(H)$ the $H$-character of the induced module $N\uparrow_K^H$. In particular, for $N=k$ viewed as $K$-module via the augmentation map $\eps_{_K}$ let $\eps\uparrow_K^H:=\eps_{ _K}\uparrow_K^H$ be the character of $H\ot_Kk=H//K$.
\begin{cor}\label{ker}
Let $K$ be a Hopf subalgebra of $H$. Then $K$ is normal in $H$ if
and only if $H_{_{\eps\uparrow_K^H}}=K$.

\end{cor}
\begin{proof}
Suppose $K$ is a normal Hopf subalgebra of $H$. With the notations
from the above theorem, since
$\eps\uparrow_K^H=|L|\pi^*(t_{_L})$ and $H_{_{|L|\pi^*(t_{ _L})}}=K$, it follows $H_{_{\eps\uparrow_K^H}}=K$.

Conversely, suppose that $H_{_{\eps\uparrow_K^H}}=K$. Then using the Remark \ref{restr} it follows
$\eps{\uparrow}_K^ H{ \downarrow }_K= \frac{|H|}{|K|}\eps_{_K}$. Using Frobenius reciprocity
this implies that for any irreducible character $\chi$ of $H$ we have that the value of
$m(\ch \dw_K, \; \eps_K)=m(\ch,\; \eps\uparrow_K^H)$ is either $\ch(1)$ if $\ch$ is a
constituent of $\eps \uw_K^H$ or $0$ otherwise. But if $\Lam_{ _K}$ is
the idempotent integral of $K$ then $m(\ch \dw_K, \;
\eps_{_K})=\ch(\Lambda_{_K})$(\cite{N}). The remark above implies that
$\Lambda_{_K}$ is a central idempotent of $H$. Then $K$ is a normal
Hopf subalgebra of $H$ by \cite{Mas'}.
\end{proof}
\section{Central characters}\label{centralch}

Let $H$ be finite dimensional semisimple Hopf algebra over $k$. Consider the central subalgebra of $H$ defined by $\hat{\mtr{Z}}(H)=Z(H) \bigcap C(H^*)$. It is the algebra of $H^*$-characters which are central in $H$.
Let $\hat{\mtr{Z}}(H^*):=Z(H^*) \bigcap C(H)$ be the dual concept, the subalgebra of  $H^*$-characters which are central in $H$.

Let $\phi:H^* \ra H$ given by $\;f \mapsto f\rightharpoondown\Lambda_{ _H} $ where $f\rightharpoondown\Lambda_{ _H}=f(S({\Lambda_{ _H}}_{ _1})){\Lambda_{ _H}}_{ _2}$. Then $\phi$ is an isomorphism of vector spaces with the inverse given by $\phi^{-1}(h)=h  \rightharpoonup |H|t_{ _H}$ for all $h \in H$ \cite{Montg}. Recall that $(h \rightharpoonup f)(a)=f(ah)$ for all $a,h\in H$ and $f \in H^*$.

\sub{} \label{formulae} With the above notations, it can be checked that $\phi(\xi_d)=\frac{\eps(d)}{|H|}d$ and $\phi^{-1}(\xi_{ _\ch})=\ch(1)\ch$ for all $d \in \mtr{Irr}(H^*)$ and $\ch \in \mtr{Irr}(H)$, (see for example \cite{Montg}). Recall that $\xi_d$ is the
central primitive idempotent of $H^*$ associated to the simple representation corresponding to the character $d$ and $\xi_{ _\ch}$ is central primitive idempotent of $H$ associated to the simple representation corresponding to the character $\ch$.

We use the following description of $\hat{\mtr{Z}}(H^*)$ and $\hat{\mtr{Z}}(H)$ which is given in \cite{Zh}.
Since $\phi(C(H))=\mtr{Z}(H)$ and $\phi(\mtr{Z}(H^*))=C(H^*)$ it follows that the restriction $$\phi| _{ _{\hat{\mtr{Z}}(H^*)}}:\hat{\mtr{Z}}(H^*)\ra\hat{\mtr{Z}}(H)$$ is an isomorphism of vector spaces .

Since $\hat{\mtr{Z}}(H^*)$ is a commutative algebra over the algebraically closed field $k$ it has a vector space basis given by its primitive idempotents. Since $\hat{\mtr{Z}}(H^*)$ is a subalgebra of $\mtr{Z}(H^*)$ each primitive idempotent of $\hat{\mtr{Z}}(H^*)$ is a sum of primitive idempotents of $\mtr{Z}(H^*)$. But the primitive idempotents of $\mtr{Z}(H^*)$ are of the form $\xi_d$ where $d \in \mtr{Irr}(H^*)$.

Thus, there is a partition $\{\mathcal{Y}_j\}_{j\in J}$ of the set of irreducible
characters of $H^*$ such that the elements $(e_j)_{j \in J}$ given by $$e_j=\sum_{d \in \mtc{Y}_j}\xi_d$$ form
a basis for $\hat{\mtr{Z}}(H^*)$. Remark that $e_j(d)=\eps(d)$ if $d \in \mtc{Y}_j$ and $e_j(d)=0$ if
$d \notin \mtc{Y}_j$.

Since $\phi(\hat{\mtr{Z}}(H^*))=\hat{\mtr{Z}}(H)$ it follows that $\hat{e_j}:=|H|\phi(e_j)$ is a basis for $\hat{\mtr{Z}}(H)$. Using first formula from \ref{formulae} one has $$\widehat{e_j}=\sum_{d \in \mtc{Y}_j}\eps(d)d$$

\sub{} \label{part} By duality, the set of irreducible
characters of $H$ can be partitioned into a finite collection of subsets
$\{\mathcal{X}_i\}_{i\in I}$ such that the elements $(f_i)_{i \in I}$ given by $$f_i=\sum_{\ch \in
\mtc{X}_i}\ch(1)\ch$$ form a $k$-basis for $\hat{\mtr{Z}}(H^*)$ and the elements $\phi(f_i)=\sum_{\ch \in
\mtc{X}_i}\xi_{ _\ch}$ form a basis for $\hat{\mtr{Z}}(H)$. Clearly $|I|=|J|$.

The following remark will be used in the proof of the next proposition.
\begin{rem} \label{polyn}Let $M$ be a representation of a  semisimple Hopf algebra $H$. Consider the set $\mathcal{C}$ of all simple representations of $H$ which are direct summands in all the tensor powers $M^{\ot\;n}$. Then $\mathcal{C}$ is closed under tensor product and $``\;^*\;"$ and it generates a Hopf algebra $L$ which is a quotient of $H$ \cite{PQ}.  Note that if $\mtc{C} \subset C(H)$ is closed under multiplication and $``\;^*\;"$ then using the dual version of \ref{hopfgen} it follows that $\mtc{C}$ generates a Hopf subalgebra $H^*_{\mtc{C}}$ of $H^*$. It follows that $L=(H^*_{\mtc{C}})^*$. If $M$ has character $\chi \in H^*$ then the character $\pi^*(t_{ _L}) \in  C(H)$ can be expressed as a polynomial in $\chi$ with rational coefficients (see Corollary 19, \cite{NR'}).
\end{rem}

\begin{prop}\label{quotient}
Suppose $\chi$ is a character of $H$ which is central in $H^*$. Then
$H_{_{\ch}}$ is a normal Hopf subalgebra of $H$ and the simple
representations of $H//H_{_{\ch}}=H/H_{_{\ch}}^+H$ are the simple constituents
of all the powers of $\ch$.
\end{prop}
\begin{proof}
Since $\ch \in \hat{Z}(H^*)$, with the above notations one has $\ch=\sum_{j=0}^s\alpha_{_j}  e_j$, where $\al_j \in k$. It follows that $\ch(d)=\alpha_{_j}\eps(d)$ if $d \in \mtc{Y}_j$. Therefore if  $d \in \mtc{Y}_j$ then $d\in \mtr{ker}\;\ch$ if and only if $\alpha_{_j}=\ch(1)$. 

This implies $\mtr{\ker}\;\ch$ is the union of all the sets
$\mtc{Y}_j$ such that $\alpha_{_j}=\ch(1)$. Using formula \ref{f2} the integral
$|H_{_{\ch}}|\Lambda_{H_{_\ch}}$ can be written as
\begin{equation*}|H_{_{\ch}}|\Lambda_{H_{_\ch}}=\sum_{d \in \mtr{ker}\;\ch}\eps(d)d=\sum_{\{j \;|\; \alpha_j=\ch(1)\}}\;\sum_{d \in \mtc{Y}_j}\eps(d)d=\sum_{\{j \;|\; \alpha_j=\ch(1)\}}\widehat{e_j}
\end{equation*}
Then $\Lambda_{H_{_{\ch}}}$ is central in $H$ since each
$\widehat{e_j}$ is central in $H$. Theorem 2.1 from
\cite{Mas'} implies that $H_{_{\ch}}$ is normal in $H$.

Let $V$ be an $H$-module with character $\ch$ and $I=\bigcap_{m \geq 0} \mtr{Ann}(V^{\ot\;m})$.  If $L$ is the quotient Hopf algebra of $H$ generated by the constituents of all the powers of $\ch$ then from \cite{KSZ} one has that $L=H/I$. Note that $I \supset HH_{_{|L|\pi^*(t_{ _L})}}^+ $.

Using Remark \ref{polyn} one has that $\pi^*(t_{_L})$ of $L$ is a polynomial in $\ch$ with rational
coefficients . Since $\chi$ is central in $H^*$ it follows that $\pi^*(t_{_L})$ is a
central element of $H^*$ and thus $L^*$ is a normal Hopf subalgebra of $H^*$.
Using \ref{normalfacts} (for $L^*  \hookrightarrow H^*$ ) it follows that $H//(H^*//L^*)^* = L$. Then if  $K=(H^*//L^*)^*$ one has $H//K=L$.
Theorem \ref{main} implies that $H_{_{|L|\pi^*(t_{ _L})}}=K$ thus $H//H_{_{|L|\pi^*(t_{ _L})}}=H//K=L$.
But $L=H/I$  and  since $I  \supset HH_{_{|L|\pi^*(t_{ _L})}}^+$ it follows that $HH_{_{|L|\pi^*(t_{ _L})}}^+=I$.

It is easy to see that $HH_{_{\ch}}^+ \subset I$ since the elements of $H_{_{\ch}}$ act as $\eps$ on each tensor power of $V$ (see \ref{prod}).

On the other hand $|L|t_{ _L}$ is the regular character of $L$. Then
$\mtr{ker}\;\ch \supset \mtr{ker}\;|L|\pi^*(t_{_L})$ since $\ch$ is a constituent of $|L|\pi^*(t_{_L})$. Thus $I \supset HH_{\ch}^+ \supset HH_{_{|L|\pi^*(t_{ _L})}}^+$.

Since $HH_{_{|L|\pi^*(t_{ _L})}}^+=I$ it follows $HH_{\ch}^+=I$ and thus $H//H_{_{\ch}}=L$.
\end{proof}

Theorem \ref{main} and the previous proposition imply the following corollary:
\begin{cor} Let $H$ be a finite dimensional semisimple Hopf. A Hopf subalgebra of $H$ is normal if and only if it is the kernel of a character $\ch$ which is central in $H^*$.
\end{cor}

Let $H_i:=H_{_{f_i}}$. From Proposition \ref{quotient} it follows that $H_i$ is a normal Hopf subalgebra of $H$. If $K$ is any other normal Hopf subalgebra of $H$ then Theorem \ref{main} implies that $K=H_{_{\ch}}$ for some central character $\ch$. Following \cite{Zh} one has $\ch=\sum_{i=1}^sm_if_i$ for some positive integers $m_i$. Then $\mtr{ker}\;\ch=\cap_{i=1}^s\mtr{ker}\;f_i$ which implies that $H_{_{\ch}}=\bigcap_{i=1}^sH_i$. Thus any normal Hopf subalgebra is an intersection of some of these Hopf algebras $H_i$.
\sub{} \label{normgen}If $K$ and $L$ are normal Hopf subalgebras of $H$ then $KL=LK$ is a normal Hopf subalgebra of $H$ that contains both $K$ and $L$. Indeed, if $a\in K,\;b\in L$  then one has $ab=(a_1bS(a_2))a_3 \in LK$ and $x_1abS(x_2)=x_1aS(x_2)x_3bS(x_4)\in KL$ for all $x \in H$.

Let $L$ be any Hopf subalgebra of $H$. We define $core(L)$ to be the
biggest Hopf subalgebra of $L$ which is normal in $H$. Based on \ref{normgen} clearly $\mtr{core}(L)$ exists and it is unique.

In the next theorem the following terminology will be used. If $M$ and $N$ are representations of the Hopf algebra $H$ affording the characters $\ch$ and $\mu$ then we say that $\ch\;is \;a\;constiutent\; of\; \mu$ if $M$ is isomorphic to a submodule of $N$.

If $A$ is a Hopf subalgebra of $H$ then there is an isomorphism of $H$-module $H/HA^+\cong H\ot_{ _A}k$ given by $\hat{h}\mapsto h\ot_{ _A}1$.
Thus if $A\subset B \subset H$ are Hopf subalgebras of $H$ then $\eps\uw_{_{B}}^H$ is a constituent of $\eps\uw_{_{A}}^H$ since there is a surjective $H$-module map $H/HA^+ \ra H/HB^+$.

\sub{}\label{centrconst} Suppose $\ch$ an $\mu$ are two characters of two representations $M$ and $N$ of $H$ such that $\mu$ is central in $H^*$ and $\ch$ is an irreducible character which is a constituent of $\mu$. Using \ref{part} it follows that $\ch \in \mathcal{X}_{i_0}$ for some $i_0 \in I$. Since $\mu$ is central in $H^*$ it follows that $\mu$ is linear combination with nonnegative integer coefficients of the elements $f_i$. Since $\ch$ is a constituent of $\mu$ it follows that $f_{i_0}$ is also a constituent of $\mu$.

\begin{rem}\label{incl} Suppose that $M$ and $N$ are two $H$-modules affording the characters $\chi$ and $\mu$. If $M$ is a submodule of $N$ then $\mtr{ker}\;\mu \subset \mtr{\ker}\;\ch$ and consequently $H_{_{\mu}}$ is a Hopf subalgebra of $H_{_{\ch}}$. Indeed, if $d \in \mtr{ker}\;\mu$ then $d$ acts as $\eps(d)\mtr{Id}_N$ on $N$. It follows that $d$ acts also as the scalar $\eps(d)$ on $M$ and thus $d \in \mtr{ker}\;\ch$.
\end{rem}

\begin{thm}If $\ch$ is an irreducible character of $H$ such that $\ch \in \mtc{X}_i$ for some $i \in I$ then $\mtr{core}(H_{_{\ch}})=H_{_{f_i}}$.
\end{thm}
\begin{proof} Since $\ch \dw_{H_{_{\ch}}}=\ch(1)\eps_{H_{\ch}}$ by Remark \ref{restr}, it follows that $\ch$ is a constituent of $\eps\uw_{H_{\ch}}^ H$. Thus $\mtr{ker}\;\ch \supseteq \mtr{ker}\;(\;\eps\uw_{H_{\ch}}^H)$ and $H_{_{\ch}} \supset H_{_{\eps\uw_{H_{\ch}}^H}}$.

Let $H_1=H_{_{\ch}}$ and
$$H_{_{s+1}}=H_{\eps\uw_{_{H_s}}^H}\;\;\;\text{for} \;\; s \geq 1.$$ The
above remark implies implies that $H_s \supseteq H_{s+1}$. Since $H$
is finite dimensional we conclude that there is $l \geq 1$ such that
$H_l=H_{l+1}=\cdots=H_{l+n}=\cdots$.

Corollary \ref{ker} gives that $H_l$
is a normal Hopf subalgebra of $H$.

We claim that
$\mtr{core}(H_{_{\ch}})=H_l$.

Indeed, for any normal Hopf
subalgebra $K$ of $H$ with $K \subset H_{_{\ch}}\subset H$ we have that $\eps\uw_{_{H_{_{\ch}}}}^H$
is a constituent of  $\eps\uw_{_K}^H$ and then using Corollary \ref{ker} it follows that
$K=H_{_{\eps\uw_{_K}^H }}\subseteq H_{_{\eps\uw_{H_{_{\ch}}}^H}}=H_2$.

Inductively, it can
be shown that $K \subset H_s$ for any $s \geq 1$, which implies that
$\mtr{core}(H_{_{\ch}})=H_l$.

It remains to show that
$H_l=H_{_{ f_i}}$.

Since $f_i$ is a central character of $H^*$ it follows from Proposition \ref{quotient} that $H_{_{f_i}}$ is a normal Hopf subalgebra of $H$. Since $\ch$ is a constituent of $f_i$ by Remark \ref{incl} one has that $H_{_{f_i}}\subset H_{_{\ch}}$. Since $H_l$ is the core of $H_{\ch}$ one gets that $H_{_{f_i}} \subset H_l$.

Next we will show that $H_l\subset H_{_{f_i}}$.
By Proposition \ref{quotient}
$t_l:=\eps\uw_{_{H_l}}^H$ is the $H$-character of $H//H_l=H/{H_l}^+H$ and 
is central in $H^*$. Since $\ch \dw_{H_{_{\ch}}}=\ch(1)\eps_{_{H_{\ch}}}$ and $H_l \subset H_{_{\ch}}$ it follows that $\ch\dw_{ _{H_l}}=\ch(1)\eps_{ _{H_l}}$. By Frobenius reciprocity one has that $\ch$ is a constituent of the character $t_l$. Using \ref{centrconst}
it follows that $f_i$ is also constituent of $t_l$. Then $\mtr{ker} \;(f_i) \supset \mtr{ker}\; (t_l)$ and $H_{_{f_i}} \supset H_{_{t_l}}=H_l$.
\end{proof}

\begin{prop}\label{bigcap}
Let $H$ be a semisimple Hopf algebra. Then
\begin{equation*}
\bigcap_{^{\ch} \in {\mtr{Irr}(H)}} \mtr{z}_{_\chi}=\bar{G}(H)
\end{equation*}
where $\bar{G}(H)$ is the set of all central grouplike elements of
$H$.
\end{prop}
\begin{proof}Any central grouplike element $g$ of $H$ acts as a scalar on each simple $H$-module. Since $g^{\mtr{exp}(H)}=1$ it follows that this scalar is a root of unity and then $$\bar{G}(H) \subset \bigcap_{\ch_{ \in Irr(H)}}
\mtr{z}_{_\chi}$$

Let $d \in \bigcap_{\ch_{ \in \mtr{Irr}(H)}}
\mtr{z}_{_\chi}$. If $C_{ _d}$ is the simple subcoalgebra of $H$ associated to $d$ (see \ref{assoc}) then $d=\sum_{i=1}^{\eps(d)}x_{ii}$.  Remark \ref{frompf} implies that
$x_{ij}$ acts as $\delta_{i,j}\al_{_\chi}Id_{M_{_\chi}}$ on
$M_{_\chi}$ where $\al_{_\chi}$ is a root of unity. For $i \neq j$,
it follows that  $x_{ij}$ acts as zero on each irreducible
representation of $H$. Therefore $x_{ij}=0$ for all $i \neq j$ and $d$ is a group like
element of $H$. Since $d$ acts as a scalar on each irreducible
representation of $H$ we have $d \in Z(H)$ and therefore $d \in
\bar{G}(H)$.

\end{proof}

\begin{rem} If $\ch \in \hat{\mtr{Z}}(H^*)$ then Proposition \ref{quotient} together with Theorem $5.3$ of \cite {KSZ} imply that $\bar{G}(H//H_{_{\ch}})$ is a cyclic group of order equal to the index of the character $\ch$.
\end{rem}
The next theorem is the generalization of the fact that
$\mtr{Z}/\mtr{ker}\;\ch$ is a cyclic subgroup of $G/\mtr{ker}\;\ch$
for any character of the finite group $G$.

\begin{thm} Let $M$ be a representation of $H$ such that its character $\chi$ is central in $H^*$.
Then $\mtr{Z}_{\ch}$ is a normal Hopf subalgebra of $H$ and
$\mtr{Z}_{_{\ch}}//H_{_{\ch}}$ is a Hopf subalgebra of
$k\bar{G}(H//H_{_{\ch}})$.
\end{thm}

\begin{proof}
Since $\ch \in \hat{\mtr{Z}}(H^*)$ one can write $\ch=\sum_{j=0}^s\alpha_j  e_j$ with $\al_j \in k$.  A similar argument to the
one in Proposition \ref{quotient} shows that
\begin{equation*}
|\mtr{Z}_{_{\ch}}|\Lambda_{_{_{\mtr{Z}_{_{\ch}}}}}= \sum_{\{j \;|\; |\alpha_j|\;=\ch(1)\}}\;\sum_{d
\in \mtc{Y}_j}\eps(d)d=\sum_{\{j \;|\; |\alpha_j|\;=\ch(1)\}}\widehat{e_j}
\end{equation*}
Therefore $\Lambda_{_{\mtr{Z}_{_{\ch}}}}$ is central in $H$ and $\mtr{Z}_{_{\ch}}$ is normal Hopf
subalgebra of $H$.  Let
$\mtr{\pi}:H \rightarrow H//H_{_{\ch}}$ be the canonical
projection. Since $H$ is free $\mtr{Z}_{ _\ch}$ module there is also an injective Hopf algebra map $i: Z_{_{\ch}}//H_{_{\ch}} \ra H//H_{_{\ch}}$ such that $i(\bar{z})=\pi(z)$ for all $z\in Z_{_{\ch}}$.

Proposition \ref{quotient} implies that all the representations of $H//H_{_{\ch}}$ are the constituents of powers of $\chi$.

From Remark \ref{prod} it follows that $\mtr{Z}_{_{\ch}} \subset
\mtr{Z}_{_{\ch^l}}$ for any nonnegative integer $l$. Then any $d \in z_{ _\ch}$ acts as a unit scalar on all the representations of $H//H_{ _\ch}$.
Thus the image of any $d\in \mtr{z}_{_{\ch}}$ under $\pi$ acts as a unit scalar on each representation of $H//H_{_{\ch}}$ and by Proposition \ref{bigcap} it is a central grouplike elements of $H//H_{_{\ch}}$. Let $C_{ _d}=<x_{ij}>$ be the coalgebra associated to $d$ as in \ref{assoc}.  Remark 1.3 also implies that $x_{ij}$ acts as zero on any tensor power of $\ch$ and then its image under $\pi$ is zero. Thus $i(Z_{ _\ch}//H_{ _\ch})\subset k\bar{G}(H//H_{ _\ch})$.

\end{proof}

\section{Double coset formula for cosemisimple Hopf algebras}\label{dcsf}

In this section let $H$ be a semisimple finite dimensional Hopf algebra as before and $K$ and $L$ be two Hopf subalgebras. Then $H$ can be decomposed as sum of $K-L$ bimodules which are free both as $K$-modules and $L$-modules and are analogues of double cosets in group theory. To the end we give an application in the situation of a unique double coset.

There is a bilinear form $m : C(H^*)\ot C(H^*) \ra k$ defined as in \cite{NR}. If $M$ and $N$ are two $H$-comodules with characters $c$ and $d$ then  $m(c,\;d)$ is defined as $\mtr{dim}_k\mtr{Hom}^H(M,\;N)$. The following properties of $m$ (see \cite{NR}) will be used later:
$$m(x, \;yz)=m(y^*,\;zx^*)=m(z^*,\;x^*y)\;\;\text{and}\;\;m(x,y)=m(y,\;x)=m(y^*,\;x^*)$$ for all $x,y,z \in C(H^*)$.

Let $H$ be a finite dimensional cosemisimple Hopf algebra and $K$,
$L$ be two Hopf subalgebras of $H$. We define an equivalence
relation on the set of simple coalgebras of $H$ as following: $C
\sim D$ if $ C \subset \mathrm{KDL}$.

Since the set of simple subcoalgebras is in bijection with $\mtr{Irr}(H^*)$ the above relation
in terms of $H^*$-characters becomes the following: $c \sim d$ if $m(c\;,\lbd_{ _K}d\lbd_{ _L}) > 0$
where $\lbd_{ _K}$ and $\lbd_{ _L}$ are the integrals of $K$ and $L$ with $\eps(\Lam_{ _K})=|K|$ and $\eps(\Lam_{ _L})=|L|$.  In order to simplify the further writing in this section we changed the notations from the previous sections.  Note that  in the previous sections $\Lam_{ _A}$ was denoting the idempotent integral of $A$ ($\eps(\Lam_{ _A})=1$) for any finite dimensional semisimple Hopf algebra $A$.

It is easy to see that $\sim$ is an equivalence
relation. Clearly $c \sim c$ for any $c \in \mtr{Irr}(H^*) $ since both $\Lam_{ _K}$ and $\Lam_{ _L}$ contain the trivial character.

Using the properties of the bilinear form $m$ given in \cite{NR}, one can see that if $c \sim d$ then $m(d,\; \Lam_{ _K}c\Lam_{ _L})=m(\Lam_{ _K}^*,\;c\Lam_{ _L}d^*)=m(c^*,\;\Lam_{ _L}d^*\Lam_{ _K})=m(c,\;\Lam_{ _K}^*d\Lam_{ _L}^*)=m(c,\;\Lam_{ _K}d\Lam_{ _L})$ since $\Lam_{ _K}^*=\Lam_{ _K}$ and $\Lam_{ _L}^*=\Lam_{ _L}$. Thus $d \sim c$.

The transitivity can be easier seen that holds in terms of simple subcoalgebras. Suppose that $c\sim d$ and $d \sim e$ and $c,\;d,$ and $e$ are three irreducible characters associated to the simple subcoalgebras $C,\;D$ and $E$ respectively. Then $C\subset KDL$ and $D \subset KEL$. The last relation implies that $KDL \subset K^2EL^2=KEL$. Thus $C \subset KEL$ and $c \sim e$.

If $\mathcal{C}_1,\mathcal{C}_2,\cdots \mathcal{C}_l$ are
the equivalence classes of $\sim$ on $\mtr{Irr}(H^*)$ then let \begin{equation}\label{pregnv}a_i=\sum\limits_{d \in
\mtc{C}_i}\eps(d)d\end{equation} for $0 \leq i \leq l$.

\begin{rem}Let $C_1$ and $C_2$ be two subcoalgebras of $H$ and $K=\sum_{n\geq 0}C_1^n$ and $L=\sum_{n \geq 0}C_2^n$ be the two Hopf subalgebras of $H$ generated by them \cite{NR}. The above equivalence relation can be written in terms of characters as follows: $c \sim d$ if
$m(c,\;c_1^ndc_2^m) > 0$ for some natural numbers $m,n\geq 0$.
\end{rem}

In the sequel, we use the Frobenius-Perron theorem for matrices with
nonnegative entries (see \cite{F}). If $A$ is such a matrix then $A$
has a positive eigenvalue $\lambda$ which has the biggest absolute
value among all the other eigenvalues of $A$. The eigenspace
corresponding to $\lambda$ has a unique vector with all entries
positive. $\lambda$ is called the principal value of $A$ and the
corresponding positive vector is called the principal vector of $A$. Also the eigenspace of $A$ corresponding to $\lam$ is called the principal eigenspace of $A$.

The following result is also needed:

\begin{prop}(\cite{F}, Proposition 5.)\label{transpose}
Let $A$ be a matrix with nonnegative entries such that $A$ and $A^t$
have the same principal eigenvalue and the same principal vector.
Then after a permutation of the rows  and the same permutation of the columns $A$ can be decomposed in diagonal blocks $A={A_1,\; A_2,\;\cdots,\;A_l}$
with each block an indecomposable matrix.
\end{prop}

For the definition of an indecomposable matrix see \cite{F}.

For any character $d \in C(H^*)$ let $L_d$ and $R_d$ be the left and
right multiplication with $d$ on $C(H^*)$.

\begin{thm}\label{mainMack}
Let $H$ be a finite dimensional semisimple Hopf algebra over the
algebraically closed field $k$ and $K$, $L$ be two Hopf subalgebras of
$H$. Consider
the linear operator $T=L_{ _{\Lam_{ _K}}}\circ R_{ _{\Lam_{ _L}}}$ on the character
ring $C(H^*)$ and $[T]$ the matrix associated to $T$ with
respect to the standard basis of $C(H^*)$ given by the
irreducible characters of $H^*$.
\begin{enumerate}
\item
The principal eigenvalue of $[T]$ is $|K||L|$.
\item
The eigenspace corresponding to the eigenvalue  $|K||L|$
has ${(a_i)}_{1 \leq i \leq l}$ as $k$- basis were $a_i$ are defined in \ref{pregnv}.
\end{enumerate}
\end{thm}

\begin{proof}
\begin{enumerate}
\item
Let $\lambda$ be the principal eigenvalue of $T$ and $v$ be the principal
eigenvector corresponding to $\lambda$. Then $\Lam_{ _K}v\Lam_{ _L}=\lambda v$.
Applying $\eps$ on both sides of this relation it follows that
$|K||L|\eps(v)=\lam \eps(v)$. But $\eps(v)>0$ since $v$ has positive entries and it follows that
$\lambda=|K||L|$.
\item
It is easy to see that the transpose of the matrix $[T]$ is  also $[T]$. To check that
let $x_1,\;\cdots,\;x_s$ be the basis of $C(H^*)$ given by the irreducible characters of $H^*$
and suppose that $T(x_i)=\sum_{j=1}^st_{ij}x_j$. Thus $t_{ij}=m(x_j,\;\Lam_{ _K}x_i\Lam_{ _L})$
and $t_{ji}=m(x_i,\;\Lam_{ _K}x_j\Lam_{ _L})=m(\Lam_{ _K}^*,\;x_j\Lam_{ _L}x_i^*)=m(x_j^
*,\;\Lam_{ _L}x_i^*\Lam_{ _K})=m(x_j, \Lam_{ _K}^*x_i\Lam_{ _L}^*)=t_{ij}$ since $\Lam_{ _K}^*=S(\Lam_{ _K})=\Lam_{ _K}$ and also $\Lam_{ _L}^*=\Lam_{ _L}$. Proposition \ref{transpose} implies that the eigenspace
corresponding to $\lambda$ is the sum of the eigenspaces of the
diagonal blocks $A_1, A_2,\cdots A_l$. Since each $A_i$ is an
indecomposable matrix it follows that the principal eigenspace of $A_i$ is one
dimensional (see \cite{F}). The rows and columns of any block $A_i$ are indexed over a subset $S_i$ of characters of $\mtr{Irr}(H^*)$. The definition of $\sim$ implies that each $S_i$ is a union of some of its equivalence classes. Since $A_i$ is indecomposable it follows that $S_i$ is just one of the equivalence classes of $\sim$.  Clearly $a_i$ is a principal eigenvector for $T$.  It follows that the eigenspace corresponding to the principal eigenvalue  $|K||L|$
has a $k$- basis given by $a_i$ with $0 \leq i \leq l$.
\end{enumerate}
\end{proof}

\begin{cor}\label{bim}
Let $H$ be a finite dimensional cosemisimple Hopf algebra and
$K,\;L$ be two Hopf subalgebras of $H$. Then $H$ can be decomposed
as
\begin{equation*}
H=\bigoplus_{i=1}^lB_i
\end{equation*}
where each $B_i$ is a $(K,L)$- bimodule free as both left $K$-module
and right $L$-module.
\end{cor}
\begin{proof}Consider as above the equivalence relation $\sim$ relative to the Hopf subalgebras $K$ and $L$.
For each equivalence class $\mathcal{C}_i$ let $B_i=\bigoplus_{C \in
\mathcal{C}_i}C$. Then $KB_iL=B_i$ from the definition of the
equivalence relation.  Then
$B_i=KB_iL \in {}_{K}\mathcal{M}_L^H$ which implies that $B_i$ is free
as left $K$-module and right $L$-module \cite{NZ}.
\end{proof}
\begin{cor}With the above notations, if $d \in \mathcal{C}_i$ then
\begin{equation}\label{formula}\frac{\Lam_{ _K}}{|K|}d\frac{\Lam_{ _L}}{|L|}=\epsilon(d)\frac{a_i}{\epsilon(a_i)}\end{equation}
\end{cor}
\begin{proof}
One has that $\Lam_{ _K}d\Lam_{ _L}$ is an eigenvector of $T=L_{ _{\Lam_{ _K}}}\circ R_{ _{\Lam_{ _L}}}$ with the maximal eigenvalue $|K||L|$. From Theorem \ref{mainMack} it follows that $\Lam_{ _K}d\Lam_{ _L}$ is a linear combination of the elements $a_j$. But $\Lam_{ _K}d\Lam_{ _L}$ cannot contain any  $a_j$ with $j \neq i$ because all the irreducible characters entering in the decomposition of the product are in $\mathcal{C}_i$. Thus $\Lam_{ _K}d\Lam_{ _L}$ is a scalar multiple of $a_i$ and the formula \ref{formula} follows.
\end{proof}

\begin{rem}\label{oneside}
Setting $C_1=k$ in Theorem \ref{mainMack} we obtain Theorem 7
\cite{NR}. The above equivalence relation becomes $c \sim d$ if and
only if $m(c,\;dc_2^m) > 0$ for some natural number $m \geq 0$. The
equivalence class corresponding to the simple coalgebra $k1$
consists of the simple subcoalgebras of the powers $C_2^m$ for $m
\geq 0$. Without loss of generality we may assume that this
equivalence class is $\mathcal{C}_1$. It follows that
\begin{equation*}
\frac{d}{\eps(d)}\frac{a_1}{\eps(a_1)}=\frac{a_i}{\eps(a_i)}
\end{equation*}for any irreducible character $d\in \mathcal{C}_i$.
\end{rem}

Let $H$ be a semisimple Hopf algebra and $A$ be a Hopf subalgebra. Define $H//A=H/HA^+$ and let $\pi : H \ra H//A$ be the module projection. Since $HA^+$ is a coideal of $H$ it follows that $H//A$ is a coalgebra and $\pi$ is also a coalgebra map.

Let $k$ be the trivial $A$-module via the counit $\eps$. It can be checked that $H//A\cong H\ot_{ _A}k$ as $H$-modules via the map $\hat{h}\ra h\ot_{ _A}1$. Thus $\mtr{dim}_kH//A=\mtr{rank}_{ _A}H$.

If $L$ and $K$ are Hopf subalgebras of $H$ define also $LK//K:=LK/LK^+$. $LK$ is a right free $K$-module since $LK \in \mtc{M}_K^H$.
A similar argument to the one above shows that $LK//K\cong LK\ot_{ _K}k$ where $k$ is the trivial $K$-module. Thus $\mtr{dim}_kLK//K=\mtr{rank}_{ _K}LK $.

\begin{thm} Let $H$ be a semisimple Hopf algebra and $K,\;L$ be two Hopf subalgebras of $H$. Then $L//L\cap K\cong LK//K$ as coalgebras and left $L$-modules.
\en{thm}

\bn{proof} Define the map $\phi: L \ra LK//K$ by $\phi(l)=\hat{l}$. Then $\phi$ is the composition of $L \hookrightarrow LK \ra LK//K$ and is a coalgebra map as well as a morphism of left $L$-modules. Moreover $\phi$ is surjective since $\widehat{lk}=\eps(k)\hat{l}$ for all $l \in L$ and $k \in K$. Clearly $L(L\cap K)^+ \subset \mtr{ker}(\phi)$ and thus $\phi$ induces a surjective map $\phi: L//L\cap K\ra LK//K$.

Next it will be shown that
\bn{equation}\label{dims}\frac{|L|}{|L\cap K|}=\frac{|LK|}{|K|}\en{equation} which implies that $\phi$ is bijective since both spaces have the same dimension.
Consider on $\mtr{Irr}(H^*)$ the  equivalence relation introduced above and corresponding to the linear operator $L_{ _{\Lam_{ _L}}}\circ R_{ _{\Lam_{ _K}}}$.
Assume without loss of generality that $\mathcal{C}_1$ is the equivalence class of the character $1$ and
put $d=1$ the trivial character, in the formula \ref{formula}. Thus $\frac{\Lam_{ _L}}{|L|}\frac{\Lam_{ _K}}{|K|}=\frac{a_1}{\epsilon(a_1)}$.  But from the definition of $\sim$ it follows that $a_1$ is formed by the characters of the coalgebra $LK$.
On the other hand $\Lam_{ _L}=\frac{1}{|L|}\sum_{d \in Irr(L^*)}\eps(d)d$ and $\Lam_{ _K}=\frac{1}{|K|}\sum_{d \in Irr(K^*)}\eps(d)d$
(see \cite{Lar}). Equality \ref{dims} follows counting the multiplicity of the
irreducible character $1$ in $\Lam_{ _K}\Lam_{ _L}$. Using \cite{NR} we know
that $m(1\;, dd')>0$ if and only if $d'=d^*$ in which case $m(1\;,
dd')=1$. Then $m(1, \;\frac{\Lam_{ _L}}{|L|}\frac{\Lam_{ _K}}{|K|})=\frac{1}{|L||K|}\sum_{d \in \mtr{Irr}(L\cap K)}\eps(d)^2=\frac{|L\cap K|}{|L||K|}$ and $m(1, \;\frac{a_1}{\epsilon(a_1)})=\frac{1}{\eps(a_1)}=\frac{1}{|LK|}$.
\end{proof}

\begin{prop}Let $H$ be a finite dimensional cosemisimple Hopf algebra and $K,\;L$ be two Hopf subalgebras of $H$.
If $M$ is a $K$-module then
\begin{equation*}
M\uw_K^{LK}\dw_L \cong (M \dw_{L\cap K})\uw^L
\end{equation*}
\end{prop}
\begin{proof}

For any
$K$-module $M$ one has
$$M\uw^{LK}\dw_L=LK\ot_KM$$
while $$(M\dw_{L\cap K})\uw^L=L\ot_{L\cap K}M$$
The previous Proposition implies that $\mtr{rank}_{_ K}LK=\mtr{rank}_{_ {L\cap K}}L$ thus both modules above have the same dimension.

Define the map $\phi: L\ot_{L\cap K}M \ra LK\ot_KM$ by $\phi(l\ot_{L\cap K} m)=l\ot_K m$ which is the composition of $L\ot_{L\cap K}M \hookrightarrow LK\ot_{L\cap K}M \rightarrow LK\ot_KM$. Clearly $\phi$ is a surjective homomorphism of $L$-modules. Equality of dimensions implies that $\phi$ is an isomorphism.

\end{proof}
If $LK=H$ then
the previous theorem is the generalization of Mackey's theorem decomposition for groups in
the situation of a unique double coset.
\bibliographystyle{amsplain}
\bibliography{n2}
\end{document}